\documentclass[12pt]{article}

\usepackage{amsmath}
\usepackage{amsfonts}
\usepackage{amssymb}
\usepackage{graphicx}
\usepackage{color}
\usepackage{natbib}
\usepackage{stackengine}
\usepackage{epsf,array,delarray,amsmath,amssymb,amsthm}
\def\tcr{}

\newcommand\barbelow[1]{\stackunder[1.2pt]{$#1$}{\rule{.8ex}{.075ex}}}

\newcommand{\half}{  {\scriptstyle \frac{1}{2} }\,  }
\newcommand{\third}{  {\scriptstyle \frac{1}{3} }\,  }

\newcommand{\R}{{\mathbb R}}
\newcommand{\sign}{  \hbox{\rm sign}  }
\newcommand{\uW}{ \barbelow{W} }
\newcommand{\uX}{ \barbelow{X} }

\newcommand{\uZ}{ \barbelow{Z} }

\newcommand{\sY}{{\mathcal Y}}
\newcommand{\sF}{ {\mathcal F} }
\newcommand{\kk}{ N }

\newtheorem{proposition}{Proposition}
\newtheorem{lemma}{Lemma}
\newtheorem{theorem}{Theorem}

\title{When is it best to follow the leader?}
\author{Philip A. Ernst\footnote{Department of Statistics, Rice University, Houston, Texas 77005, U.S.A.}\,, L.C.G. Rogers\footnote{Statistical Laboratory, Center for Mathematical Sciences, Cambridge CBS 0WB, U.K.}\,, and Quan Zhou\footnote{Department of Statistics, Rice University, Houston, Texas 77005, U.S.A.}}
\date{\today}
\begin{document}
\maketitle

\begin{abstract}
An object is hidden in one of $N$ boxes. Initially, the probability that it is in box $i$ is $\pi_i(0)$. You then search in continuous time, observing box $J_t$ at time $t$, and receiving a signal as you observe: if the box you are observing does not contain the object, your signal is a Brownian motion, but if it does contain the object your signal is a Brownian motion with positive drift $\mu$. It is straightforward to derive the evolution of the posterior distribution $\pi(t)$ for the location of the object. If $T$ denotes the first time that one of the $\pi_j(t)$ reaches a desired threshold $1-\varepsilon$, then the goal is to find a search policy $(J_t)_{t \geq 0}$ which minimizes the mean of $T$. This problem was studied by  \citet{posner1966continuous} and by \cite{zigangirov1966problem}, who derive an expression for the mean time of a conjectured optimal policy, which we call {\em follow the leader} (FTL); at all times, observe the box with the highest posterior probability. Posner \& Rumsey assert without proof that this is optimal, and Zigangirov offers a proof that if the prior distribution is uniform then FTL is optimal. In this paper, we show that if the prior is not uniform, then FTL is {\em not} always optimal; for uniform prior, the question remains open.\\

\end{abstract}

\section{Introduction.}\label{intro}
This paper studies a classical search problem first considered independently by \cite{posner1966continuous} and \cite{zigangirov1966problem}. An object is hidden in one of $N$ boxes; we denote by $j^*$ the index of the true box. Initially,
\[
P( j^* = i ) = \pi_i(0).
\]
We then observe in continuous time, choosing to search box $J_t$ at time $t$. We see a signal process $Y$ whose dynamics are
\begin{equation}
dY(t) = dW(t) + \mu I_{ \{ J_t = j^*\} } \, dt,
\label{dY}
\end{equation}
where $\mu>0$ is a known constant, and $W$ is a Brownian motion. It is straightforward to derive the evolution of $\pi(t)$, the posterior distribution at time $t$. The objective proposed by \cite{posner1966continuous} is to choose $(J_t)_{t \geq 0}$ to minimize $ET$, where 
\begin{equation}
T = \inf\{t:  \max_j \pi_j(t) \geq 1- \varepsilon\},
\label{Tdef}
\end{equation}
where $\varepsilon \in (\half,1)$ is some desired error bound.

To the best of our knowledge, the solution to this specific problem was studied by three papers: \cite{posner1966continuous, zigangirov1966problem} and\\~\cite{klimko1975optimal}. 
\cite{posner1966continuous} asserted, without proof, that the optimal strategy is to always search the box with the largest posterior probability. We call this policy the \textit{follow the leader} (FTL) policy.
They formulated the FTL strategy as the limit of a sequence of discrete-time approximations, which was later shown by~\cite{klimko1971some} not to be tight.
\cite{zigangirov1966problem} considered only the case of uniform prior distribution, that is,   $\pi_i(0) = 1/N$  for $i=1,\dots, N$, and offered a proof for the optimality of FTL. 
However, this proof lacks clarity on a number of points, and we were not able to verify the arguments given.
\cite{klimko1975optimal} provided a proof for the optimality of FTL for arbitrary prior distribution, but, as we will explain later, their proof is in error.
The main result of our paper is to give counterexamples that clearly show FTL is not optimal for some specific values of  $(\pi_1(0), \dots, \pi_N(0))$. An additional contribution  is the characterization of the solution to a class of stochastic differential equations, which plays a key role in our calculations, and can be considered to be generalizations of Tanaka's SDE. 

\subsection{Literature review.}
Optimal scanning problems apparently date back to~\citet{shiryaev1964theory}. Our specific problem of interest was considered in the works of \cite{posner1966continuous, zigangirov1966problem} and~\cite{klimko1975optimal}. We now briefly review other variants of the problem that are closely related to this work. \\
\indent \cite{liptser1965bayesian} considered a setup with two boxes ($N=2$) and allow for the possibility that the object may not be in either box. The task is to determine if the missing object is in one of the two boxes. 
\cite{dragalin1996simple} considered general stochastic processes other than Brownian motion and proposed a scanning rule based on the sequential probability ratio test of~\cite{wald1945sequential}. \\
\indent Another class of problems similar to optimal scanning is problems of ``quickest search.'' 
These problems are often formulated under the setting $N\rightarrow \infty$ with an unknown number of boxes containing the hidden objects. We refer readers to~\cite{lai2011quickest} for a discrete-time solution and \cite{bayraktar2014quickest} for a continuous-time solution to these quickest search problems.  We also note that optimal scanning problems and quickest search problems are known collectively as ``screening problems.''
We refer readers to the references given in~\cite{heydari2016quickest} for other variants of screening problems. 
More generally, such problems can be viewed as sequential decision problems. 
References on this topic include~\cite{dvoretzky1953sequential}, \citet[Chap. 4]{shiryaev2007optimal}, and \citet[Chap. VI]{peskir2006optimal}.
 

\section{The evolution of the posterior.} \label{S1}
As we noted at \eqref{dY}, the signal process $Y$ evolves as
\begin{equation*}
dY(t) = dW(t) + \mu I_{ \{ J_t = j^* \} } \, dt.
\end{equation*}
If $(\sY_t)_{t \geq 0}$ is the filtration of the observation process, 
and we choose to search box $J_t$ at time $t$, then the posterior likelihood (relative to Wiener measure) that the true box is $j$, given $\sY_t$, is  
\begin{equation}
z_j(t) = \pi_j(0) \; \exp\biggl(\;
\int_0^t \mu I_{\{J_s = j\}} \; dY_s - \half \mu^2 \int_0^t I_{\{J_s = j\}} \; ds
\;\biggr).
\label{zjdef}
\end{equation}
The posterior probabilities are obtained by normalizing the $z_j$:
\begin{equation}
\pi_j(t) = z_j(t) / \bar{z}(t),
\label{pi_t}
\end{equation}
where of course $\bar{z}(t) = \sum_j z_j(t)$.  Now
\begin{equation}
dz_j(t)  = z_j(t) \, \mu I_{\{J_t = j\}} \, dY_t
\label{dzj}
\end{equation}
so if we write $X_j(t) = \mu^{-1} \log( z_j(t) )$ then we have
\begin{equation}
dX_j(t) = I_{\{J_t = j\}} \,( dY_t - \half \mu dt).
\label{dxj}
\end{equation}
The evolution \eqref{dY} of $Y$ is expressed in the filtration of $W$, but familiar results of filtering theory (see \cite{kallianpur1972stochastic}) establish that we can rewrite the evolution in the filtration $\sY$ as
\begin{equation}
dY(t) = d\hat{W}(t) + \mu \,\pi_{J_t}(t) \, dt
\label{dY2}
\end{equation}
where $\hat{W}$ is the innovations process, a $\sY$-Brownian motion, and $\mu \, \pi_{J_t}(t)$ is the $\sY$-optional projection of the drift $\mu I_{ \{J_t = j^*\}}$ of \eqref{dY}.

\medskip
Formulating the dynamics slightly more generally, as
\begin{equation}
dz_j(t) = z_j(t) \theta_j(t) \, dY_t
\label{dzj_bis}
\end{equation}
where $\theta$ is a bounded previsible $N$-vector process, we can consider the evolution of $\pi(t)$ defined in terms of $z(t)$ by \eqref{pi_t}. Some routine calculations with It\^o's formula give us
\begin{equation}
d\pi_j(t)  = \pi_j(t) \{\;  \theta_j(t) - \theta(t) \cdot \pi(t) \;\} \{\; dY_t - \theta(t) \cdot \pi(t) \, dt \;\}, 
\label{dpi_t}
\end{equation}
\tcr{where $\theta(t) \cdot \pi(t) = \sum_i \theta_i(t) \pi_i(t)$. }
In the case of special interest to us, where $\theta_j(t)  = \mu I_{\{J_t = j\}}$, the representation \eqref{dY2} combines with \eqref{dpi_t} to show that
\begin{equation}
d\pi_j(t)  = \pi_j(t) \{\;  \theta_j(t) - \theta(t) \cdot \pi(t) \;\} d\hat{W}(t).
\label{dpi_t_bis}
\end{equation}
In particular, $\pi(t)$ is a $\sY$-local martingale; but we know this already, because $\pi_j(t)  = P (\; j^* = j  \;| \; \sY_t\,)$, which is even a martingale.

\section{The FTL policy.}\label{S2}
For the FTL policy, the dynamics \eqref{dzj_bis} has the special form
\begin{equation}
 \theta_j(t) = \mu I_j ( X(t)),
 \label{theta_FTL}
\end{equation}
where we define for $x \in \R^N$
\begin{eqnarray}
I_j(x) &=& 1 \qquad \hbox{\rm if $x_j = \max\{ x_i :i=1, \ldots,j\} > 
\max\{ x_i: i >j\} $}
\nonumber
\\
&=& 0\qquad \hbox{\rm else.}
\label{Idef}
\end{eqnarray}
Thus $I_j(x)$ is the index where the $N$-vector $x$ is maximal, taking care to avoid ambiguities when there are ties, and to ensure that $\sum_j I_j(x) = 1$. 

In these terms, the evolution of   $X_j(t)$ can be expressed as \footnote{In \cite{posner1966continuous}, the signal has a constant volatility $\sigma$,  as well as the drift $\mu$. We could always scale the signal to turn the volatility to 1, and indeed we could also replace $\mu$ by any desired positive value; this is equivalent to a constant rescaling of time, which will not affect optimality of a search policy. }
\begin{eqnarray}
dX_j(t) &=& I_j(X(t)) \, (dY_t - \half \mu \, dt)
\nonumber
\\
&=& I_j(X(t)) \, ( d\hat{W}(t) + \mu( \pi_j(t) - \half) \, dt )
\label{SDE1}
\\
&\equiv & I_j(X(t)) \, dZ_t,
\label{Zdef}
\end{eqnarray}
say. Since 
\begin{equation}
\pi_j(t) = \frac{e^{\mu X_j(t)} }{\sum_{i=1}^N e^{\mu X_i(t)}}
\label{pi_X}
\end{equation}
is a function of $X(t)$, the SDE \eqref{SDE1} is an autonomous SDE, but the coefficients are not Lipschitz, or even continuous, so the sense in which the SDE has a solution needs to be clarified.

\bigbreak

 We shall address this by firstly studying the SDE \eqref{SDE1} without the drift term:
\begin{equation}
dX_j(t)  = I_j(X(t))  \, dW_t,
\label{SDE0}
\end{equation}
where $W$ is a standard Brownian motion.  
To appreciate the issues involved, let us first consider the case $N=2$, when the SDE is 
\begin{eqnarray}
dX_1(t) &=& I_{ \{ X_1(t) > X_2(t) \} } \, dW_t = I_{ \{ X_1(t)-X_2(t) >0 \} } \, dW_t
\label{dX1}
\\
dX_2(t) &=& I_{ \{ X_2(t) \geq X_1(t)\} } , dW_t= I_{ \{ X_1(t)-X_2(t) \leq 0 \} } \, dW_t.
\label{dX2}
\end{eqnarray}
So if $Y_t \equiv X_1(t) - X_2(t)$ we have the  celebrated Tanaka SDE
\begin{equation}
dY_t = \sign(Y_t) \, dW_t
\label{dY_tanaka}
\end{equation}
for $Y$, where the definition of $\sign $ is the correct one for the definition of semimartingale local time - see Theorem IV.43.3 of \cite{rogers2000diffusions}. There is no strong solution to this SDE, but there is a weak solution, represented by taking $Y$ to be a Brownian motion started at $y_0 = X_1(0) - X_2(0)$, which we may as well suppose is positive, and then defining
\begin{equation}
 dW_t = \sign(Y_t) \, dY_t = d|Y_t| - dL_t
 \label{dWdef}
\end{equation}
where $L$ is the local time of $Y$ at zero. Then we have
\begin{equation}
dX_1(t) = I_{\{Y_t>0\}} \, dY_t = d(Y_t^+) - \half dL_t.
\end{equation}
Thus
\begin{eqnarray}
X_1(t) &=& X_2(0) + Y_t^+ - \half L_t
\label{X1fromY}
\\
X_2(t) &=& X_2(0) + Y_t^- - \half L_t.
\label{X2fromY}
\end{eqnarray}
In view of the above, we realize:
\begin{itemize}
\item We cannot hope for \eqref{SDE0} to have a strong solution;
\item We might obtain uniqueness in law for all initial values;
\item If all the $X_j$ start from 0, the sum $\sum_j X_j(t)$ is a Brownian motion when $N=2$;
\item If all the $X_j$ start from 0, the running minimum $\uX_j(t) \equiv \inf\{X_j(s): s \leq t\}$ is the same for all $j$ when $N=2$.
\end{itemize}
For general $N$, we have the analogous conclusions.

\begin{theorem}\label{thm1}
For all starting values $X(0)$, the SDE \eqref{SDE0} has a weak solution which is unique in law. If $X(0) = 0$, then

\begin{itemize}
\item[(1)] $\sum_{j=1}^N X_j(t) \equiv W_t$ is a Brownian motion;
\item[(2)] the running minimum processes coincide:
\begin{equation}
\uX_j(t) \equiv \inf\{X_j(s): s \leq t\} = N^{-1} \uW(t) \equiv \uX(t)
\end{equation}
\item[(3)] for all $t \geq 0$, 
\begin{equation}
X_j(t) = \uX(t) \qquad\hbox{\it for all but at most one index $j$};
\end{equation}
\item[(4)] $X_j(t) - \uX(t) = W(t) - \uW(t)$ if $X_j(t) > \uX(t)$;
\item[(5)] The process
\begin{equation}
X_k(t) - (N-1)^{-1} \sum_{j \neq k} X_j(t)
= \frac{NI_k(X(t)) -1}{N-1}  (W(t) - \uW(t) )
\end{equation}
is a martingale.
\end{itemize}
\end{theorem}

Theorem \ref{thm1} is proved in Appendix \ref{app}. It deals immediately with the question of existence and uniqueness of solutions of \eqref{SDE1}, because any weak solution to \eqref{SDE0} can be transformed by change of measure into a weak solution to \eqref{SDE1}, and {\it vice versa}.

\section{The value of FTL.}\label{S3}
Our aim in this section is to discover the value function  $V(x_1, \ldots, x_N) $  of the FTL policy, where    $x_j = \mu^{-1} \log z_j(0)$ denotes the initial value of the process $X_j(t)$. 
\tcr{Formally,
$$V(x_1, \ldots, x_N) = E \left( T  \mid X_1(0) = x_1, \dots, X_N(0) = x_N,  \text{policy = FTL} \right),$$ 
where $T$ is as in~\eqref{Tdef}.}
It does not appear possible to express this in closed form, but we can find a recursive algorithm for computing the value numerically. It is obvious that $V$ is a symmetric function of its arguments, so we will make the convention in what follows that the arguments of $V$ have been arranged in decreasing order:
\begin{equation}
x_1 \geq x_2 \geq \ldots \geq x_N.
\label{x_decr}
\end{equation}
It is also obvious that the value will not be changed if we add the same constant to all the arguments.

\medskip

Now suppose that all the inequalities in \eqref{x_decr} are strict, and we apply the FTL rule. What happens is that initially we observe the most likely box, box 1, up until the time $\tau_2$ when $X_1$ first falls\footnote{Of course, the termination condition may have been achieved before that time.} to $x_2$. At that time, $X_2$ begins to move, and  in accordance with Theorem \ref{thm1} we find that 
\begin{equation}
d\uX_1(t) = d\uX_2(t) = \half d\uZ_t, \,\, \max\{X_1(t), X_2(t)\}
-\uX_1(t) = Z_t - \uZ_t.
\end{equation}
This continues until the first time $\tau_3$ that one of $X_1$, $X_2$ (and hence both) falls to $x_3$. Thereafter we observe the boxes 1, 2, 3 with
\small
\begin{equation*}
d\uX_1(t) = d\uX_2(t) = d\uX_3(t) = \third d\uZ_t, \qquad \max\{X_1(t), X_2(t), X_3(t) \}
-\uX_1(t) = Z_t - \uZ_t.
\end{equation*}
\normalsize
This continues sequentially, with  the  $X_j(t)$ starting to move one after another, until all the  $X_j(t)$  achieve a common minimum $x_N$ at time $\tau_N$, or, of course, the termination criterion is achieved. So we see that calculation  of the value comes down to solving a sequence of first-exit problems, which we formalize in the following result.

\begin{theorem}\label{thm2}
Suppose that $x_1> x_2 > \ldots>x_N$, and let $V_n(x, x_{n+1}, \ldots, x_N)$ denote the value if we start with $x_1 = x_2 = \ldots = x_n = x > x_{n+1} > \ldots > x_N$. Then the values $V_n(x_n,x_{n+1}, \ldots, x_N)$ can be calculated recursively as
\begin{equation}
V_n(x_n,x_{n+1}, \ldots, x_N) = A_n(x_n) + \frac{B_n(x_n)}{1+ K_n(x_n)},
\label{Vn0}
\end{equation}
where
\begin{equation}
K_n(y) \equiv n-1+b_ne^{-\mu y}, \qquad  b_n \equiv e^{\mu x_{n+1}} + \ldots + e^{\mu x_N},
\label{Kndef}
\end{equation}
where $B_n$ is the solution to the ODE ($p_0 \equiv 1-\varepsilon$)
\small
\begin{eqnarray}
(1-p_0(1+K_n(y)) )B'(y)
&=& \tcr{ (n-1)\mu B(y) + \frac{2n(K_n(y)-1)}{\mu}   } 
\nonumber \\
&& - \frac{2 (1-2p_0) }{\mu K_n(y)}
(K_n(y)-n+1)(K_n(y)+1),
\label{B_ode}
\end{eqnarray}
\normalsize
with boundary condition
\begin{equation}
B(x_{n+1}) = \frac{(1+K_n(x_{n+1}))\{\mu V_{n+1}(x_{n+1}, \ldots,x_N)-
2 q_n(x_{n+1})(1-2\varepsilon)\} }{\mu\{ 1-p_0(1+K_n(x_{n+1})  )\} },
\label{B_bc}
\end{equation}
and where
\begin{eqnarray}
0 &=&   A_n(y) + (1-\varepsilon)B_n(y) -2q_n(y)(1-2\varepsilon)/\mu  
\label{AfromB}
\\
q_n(y) &\equiv & \mu^{-1} \log\{ \; (1-\varepsilon)(n-1+b_ne^{-\mu y})/\varepsilon \; \}.
\end{eqnarray}
The recursion begins at $n=N$ with the Posner-Rumsey value:
\begin{equation}
V_N(x_N) = M_{PR} \equiv \frac{2}{\mu^2} \biggl[ \frac{N-2}{N-1}\;(Np_0-1)
+(2p_0-1) \log\biggl(\frac{1-\varepsilon}{\varepsilon/(N-1)}
\biggr)\biggr].
\label{MPR}
\end{equation}

\end{theorem}

\begin{proof}
 Holding $x_{n+1}, \ldots, x_N$ fixed, we let $f(s,y)$ denote the value when the running minimum of $X_1, \ldots,X_n$ is $y \in [x_{n+1},x_n]$, and the unique leading particle is at $y+s \geq y$. The success criterion will be satisfied if
\begin{equation}
\frac{e^{\mu(y+s)}}{ e^{\mu(y+s)} + (n-1) e^{\mu y} + b_n}  = 1-\varepsilon,
\label{succ1}
\end{equation}
where $b_n$ is given by \eqref{Kndef}.
Hence successful termination occurs when 
\begin{equation}
s = \mu^{-1} \log\{ \; (1-\varepsilon)(n-1+b_ne^{-\mu y})/\varepsilon \; \}
\equiv q_n(y).
\label{succ2}
\end{equation}
Now up until $\tau_{n+1}$, in terms of $Z$ we have that
\begin{equation*}
f(Z_t - \uZ_t, n^{-1}\uZ_t) + t \quad\hbox{\rm is a martingale.}
\end{equation*}
The probability that we are viewing the correct box is
\begin{eqnarray}
p(s,y) &=& \frac{e^{\mu(y+s)}}{ e^{\mu(y+s)} + (n-1) e^{\mu y} + b_n}
\nonumber\\
&=& \frac{e^{\mu s}}{ e^{\mu s} + (n-1)  + b_n e^{-\mu y}}
\nonumber \\
&\equiv & \frac{e^{\mu s}}{ e^{\mu s} + K_n(y)} \; ,
\label{pdef}
\end{eqnarray}
where $K_n$ is defined at \eqref{Kndef}.
Using It\^o's formula, we arrive at the equations
\begin{eqnarray}
0 &=& \half f_{ss} + (p(s,y)-\half) \mu f_s + 1,
\label{pde1}
\\
0 &=& -  n f_s + f_y \qquad (s=0)
\label{pde2}
\end{eqnarray}
with boundary conditions
\begin{eqnarray}
f(q_n(y),y) &=& 0 \qquad\forall y \in [x_{n+1}, x_n],
\label{bc1}
\\
f(0, x_{n+1}) &=& V_{n+1}(x_{n+1}, \ldots, x_N).
\label{bc2}
\end{eqnarray}
We see that \eqref{pde1} is a second-order linear ODE in the variable $s$, whose general solution can be shown by routine calculations to be
\begin{equation}
f(s,y) = \tilde{A}(y) + \frac{\tilde{B}(y)}{e^{\mu s} + K_n(y)} + \frac{2s}{\mu} \; (1-2 p(s,y)).
\label{f_form}
\end{equation}
for some functions $\tilde{A}$, $\tilde{B}$. Equivalently, we may express the solution as
\begin{equation}
f(s,y) = A(y) + B(y) p(s,y) + \frac{2s}{\mu} \; (1-2 p(s,y))
\label{f_form2}
\end{equation}
for some functions $A$, $B$.
 The boundary condition at reflection \eqref{pde2} leads to the equation
\begin{equation}
\frac{(n-1)\mu B(y)}{1+K_n(y)}+ \frac{2n(K_n(y)-1)}{\mu(K_n(y)+1)} = A'(y) + \frac{B'(y)}{1+K_n(y)}
\label{de1}
\end{equation}
and the boundary condition \eqref{bc1} gives us
\begin{equation}
0 = A(y) + (1-\varepsilon)B(y) -2q_n(y)(1-2\varepsilon)/\mu.
\label{AB}
\end{equation}
This allows us to express $A(y)$ as a function of $y$ and $B(y)$, reducing the ODE \eqref{de1} to a first-order linear ODE for $B$. From \eqref{AB} we find that
\begin{equation}
A'(y) = -p_0 B'(y) + \frac{2 (1-2p_0) }{\mu K_n(y)} \; (K_n(y)-n+1).
\label{AB2}
\end{equation}
Returning this to \eqref{de1} leads to the first-order ODE  \eqref{B_ode} for $B$ in $y \geq x_{n+1}$.
The boundary condition \eqref{bc2} together with \eqref{AB} becomes the boundary condition \eqref{B_bc}.

\bigbreak
At the final stage, $V_N$ is a function of just one variable, and $p(s,y)$ is independent of $y$; the form \eqref{f_form2} collapses to 
\begin{equation}
f(s) = A_N + B_N\,p(s,0) + \frac{2s}{\mu} (1-2p(s,0))
\end{equation}
with the boundary conditions
\begin{equation}
 f(q_N(0) ) = 0, \qquad f'(0) = 0.
\end{equation}
Solving this for $A_N$, $B_N$ leads to the Posner-Rumsey solution \eqref{MPR}. 
\tcr{Finally, the expression~\eqref{Vn0} is obtained by letting $s = 0$.}

\end{proof}


\section{Counterexamples for the optimality of FTL.}\label{sec:counter}
In this section, we introduce an alternative strategy that can beat FTL in some circumstances.
We make use of a  classical result for the exit time of a Brownian motion at two boundaries. 
It was apparently first derived by~\cite{darling1953first}.

\begin{lemma}[\citet{darling1953first}]\label{lm:exit}
Let $W(t)$ be a Brownian motion with drift $\lambda t$ and variance $\sigma^2 t$, and started at $x$. 
Let $\rho = \lambda / \sigma^2$. 
Consider the boundaries $a$ and $b$ such that $a > x > b$. 
Then exit at one of the boundaries occurs with probability $1$, and the probability of exit at $a$ is given by 
\begin{align*}
P(x, a, b, \lambda, \sigma^2) =  \dfrac{ e^{-2\rho b} - e^{-2\rho x} }{ e^{-2\rho b} - e^{-2\rho a}  }. 
\end{align*}
Conditional on exiting at $a$, the expected time is given by 
\begin{align*}
F_a(x, a, b, \lambda, \sigma^2) = \dfrac{1}{\lambda} \left[ (a - x) + \dfrac{2(a-b) e^{- 2\rho a }  }{ e^{- 2 \rho b}  - e^{- 2\rho a } } - 
\dfrac{2(x-b) e^{- 2\rho x }  }{ e^{- 2 \rho b}  - e^{- 2\rho x } } \right]. 
\end{align*}
Conditional on exiting at $b$, the expected time is given by 
\begin{align*}
F_b(x, a, b, \lambda, \sigma^2) = \dfrac{1}{\lambda} \left[ (x  - b) + \dfrac{2(a-b) e^{- 2\rho a }  }{ e^{- 2 \rho b}  - e^{- 2\rho a } } - 
\dfrac{2(a - x) e^{- 2\rho a}  }{ e^{- 2 \rho x}  - e^{- 2\rho a } } \right] . 
\end{align*} 
\end{lemma}

\medskip\noindent
{\sc Remark.}
Observe that $F_a(x, a, b, \lambda, \sigma^2) =  F_a(x, a, b, -\lambda, \sigma^2)$. 
The same equality holds true for function $F_b$. 
This would be useful since, by~\eqref{dxj}, the drift of 
\tcr{$X_{J_t}(t)$} is either $\mu / 2$ or $-\mu/2$.

\medskip

Now we introduce an \textit{alternative strategy}, which we call ``Strategy B" .
For simplicity, let us consider three boxes with initial values $x_1 > x_2 > x_3$ (and thus prior probabilities 
 $\pi_1(0) > \pi_2(0) > \pi_3(0)$).  We shall suppose that
\begin{equation}
\pi_1(0) \equiv \frac{e^{\mu x_1}}{e^{\mu x_1} + e^{\mu x_2} +e^{\mu x_3} }
< 1-\varepsilon < \frac{e^{\mu x_1}}{e^{\mu x_1}  +2e^{\mu x_3} },
\label{inequalities}
\end{equation}
so that there exists a unique $ a \in (x_3, x_2)$ such that
\begin{equation}
\frac{e^{\mu x_1}}{e^{\mu x_1} + e^{\mu a} +e^{\mu x_3} } = 1-\varepsilon.
\label{a_def}
\end{equation}
Strategy B observes $X_2(t)$ until it reaches $a$ or $x_1$. If $X_2$ reaches $a$ before $x_1$, then the objective is achieved, in view of \eqref{a_def}. Because of Lemma \ref{lm:exit}, we know  the mean of this stopping time, and the probability that exit happens at $x_1$.    Otherwise, if $X_2$ reaches $x_1$ before $a$, we now continue with the FTL policy, whose mean remaining time to finish will be $V(x_1, x_1, x_3)$, which can be calculated according to Theorem \ref{thm2}.

By fixing $\mu = 1$, $x_3 = 0$ and searching over the corresponding domain of $(x_1, x_2)$, we obtain a few counterexamples for different values of  $\varepsilon$. These are presented in Table~\ref{table:counter}.

\begin{table}[h!]
\begin{center}
\begin{tabular}{cccccc}
\hline 
$\varepsilon$  &  $x_1$  & $x_2$      & $E_A(T) \times 10^{2}  $  &  $E_B(T) \times 10^{2} $  \\
\hline 
0.4  &  2 &  1.4  &   $3.633  $ &  $3.464 $ \\
0.3  &  2.7 & 1.7  &    $3.053 $ &  $2.936  $  \\ 
0.2  &  4.05 &  2.6 &    $1.832  $ &  $1.797  $ \\ 
0.1  & 6.2   &  4.0  &    $3.749  $ &  $3.738  $ \\
0.05  &  10.3  & 7.4  &    $10.6482 $ &  $10.6476 $   \\   
\hline 
\end{tabular} 
\caption{Counterexamples for the optimality of FTL with $N = 3$, $\mu = 1$ and $x_1 > x_2 > x_3 = 0$. 
$E_A(T)$ denotes the expected search time of the FTL strategy, which can be computed using Theorem \ref{thm2}. 
$E_B(T)$ denotes the expected search time of Strategy B, which can be computed as
explained above.}\label{table:counter}
\end{center} 
\end{table}


\section{Further discussion.}\label{sec:disc}

\subsection{Discussion of the work of~\cite{klimko1975optimal}. }\label{sec:ky}
\cite{klimko1975optimal} gave a proof for the optimality of FTL for arbitrary prior distribution; 
however, according to our counterexamples, this conclusion cannot be correct. 
Here we explain why their proof is incorrect.

Consider that we start the search by choosing one box to observe until $\pi_1(t)$ reaches either $\pi_1(0) + \alpha$ or $\pi_1(0) - \beta$, assuming $\pi_1(0) > \cdots > \pi_N(0)$. 
Denote this stopping time by $\tau$ and assume that there is no switch of the observed box before $\tau$.  
Thus for $t \in [0, \tau]$, there are $N$ possible search rules. 
In both the proofs of Lemma 3.4 and Theorem 3.5 in \cite{klimko1975optimal}, the authors assume that the posterior distribution at $\tau$ is \textit{independent of the search rule}, which is incorrect.
For example, their proof for case (i) of Theorem 3.5 relies on this incorrect assumption (the authors write that ``furthermore, the posterior distribution at exit time are also independent of the rule used.'')

For a concrete example, consider $N = 3$ and $\pi(0) = (0.5, 0.25, 0.25)$.
We wish to use this example to further explain why the reasoning of \cite{klimko1975optimal} is incorrect.
Let $\tau_1$ be the exit time of $\pi_1$ at either $1 - \varepsilon$ or $0.4$, and assume no switch of observed box before $\tau_1$. 
Lemma 3.1 of \cite{klimko1975optimal} correctly states that in order to minimize $E (\tau_1)$ we should choose to observe box $1$.
However, if box $1$ is observed and $\pi_1$ exits at $0.4$, at $\tau_1$ we have the posterior probability  $\pi(\tau_1) = (0.4, 0.3, 0.3) \equiv \pi_A$; if box $2$ is observed and and $\pi_1$ exits at $0.4$,  we have  $\pi(\tau_1) = (0.4, 0.4, 0.2) \equiv \pi_B$. 
The inductive argument used in the proof of Lemma 3.4 of  \cite{klimko1975optimal} now fails because it is no longer clear whether $\pi_A$ or $\pi_B$ would lead to a smaller expected search time after $\tau_1$. 
In fact, according to our numerics,  $\pi_B$ would give a smaller expected search time if FTL is applied.

\subsection{Open problems.}\label{sec:equal.prior}
Our work gives rise to several open problems. 
First, what is the optimal strategy for this optimal scanning problem for any prior distribution? 
For decades, it has been (incorrectly) assumed that FTL is optimal. Indeed, as we have shown, FTL is sub-optimal at least for some values of $(\pi_1(0), \dots, \pi_N(0))$. 

Another open problem concerns whether FTL is optimal for the case of uniform prior distribution. 
There is already a proof given by \cite{zigangirov1966problem}, but as already mentioned, we found the argument presented to lack clarity in various places, so we cannot be confident that the result is established.

So the answer to the question in the title is, `We don't know, but we know that it is not {\em always} best to follow the leader!'
\newpage

\section{References}
\bibliographystyle{plainnat}
\bibliography{ref}

\pagebreak
\appendix
\section{The SDE $dX = I(X) \; dW.$}\label{app}
The aim in this appendix is to prove Theorem 1. Until further notice, we shall focus on the case $X_j(0) = 0$ for all $j$. To begin with, we present  two results about any solution of the SDE \eqref{SDE0}
\begin{equation}
dX_j(t) = I_j(X(t))\, dW_t, \qquad X_j(0) = 0 \quad \forall j.
\label{SDE0_bis}
\end{equation}

\begin{proposition}\label{prop1}
For all $i$, $\uX_t^i = \uX^1_t \equiv \uX_t$.
\end{proposition}
\begin{proof}
For any $a<0$, we let $H_j(a) = \inf \{t: X_j(t) \leq a \}$. If it were the case that for some $i \neq j$ we have $X_i(H_j(a)) > X_j(H_j(a))=a$, then there has to be some time interval $(s,u)$ containing $H_j(a)$ throughout which $X_i(t)>a$. This means that throughout $(s,u)$ the process $X_j$ is not the leader, so it does not move. This contradicts the definition of $H_j(a)$. Therefore $X_i(H_j(a)) \leq a$ for all $i \neq j$, and hence $H_i(a) \leq H_j(a)$ for all $i \neq j$. Since we can interchange the r\^oles of $i$ and $j$, it must be that $H_i(a) = H_j(a)$ for all $i \neq j$, and the result follows.

\end{proof}

\begin{proposition}\label{prop2}
For all $t \geq 0$, $\uW_t = \kk \uX_t$.
\end{proposition}
\begin{proof}
Observe that $\sum_{j=1}^N X_j(t) = W_t$, since $\sum_j  I_j(x) \equiv 1$, which {\em proves statement (1) of Theorem \ref{thm1}.}
With the notation of the proof of Proposition \ref{prop1}, we have that for any $a<0$
\[
X_j( H(a) ) = \uX(H(a)) =  a \quad\forall j,
\]
where $H(a) $ denotes the common value $H_j(a)$.  Therefore $W(H(a)) = Na$. Further, for any $t < H(a)$ we have $X_j(t) > a$, and thus $W(t)>Na$. So it must be that $H(a) = \inf \{t: W(t) \leq Na \}$, and the result follows.

\end{proof}

\medskip\noindent{\sc Remark. } {\em Statement (2)  of Theorem \ref{thm1} is now proved.}

\begin{proposition}\label{prop3}
On a suitable probability space, a solution to \eqref{SDE0_bis} may be constructed.
\end{proposition}
\begin{proof}
Consider the It\^o excursion point process description of Brownian motion, using notation and terminology from Ch VI of \cite{rogers2000diffusions}. According to Proposition VI.51.2, the rate of excursions which get at least $x>0$ away from zero is
\begin{equation}
n(\{ f \in U : \sup_t |f(t)| > x\}) = 1/x,
\end{equation}
and the full excursion law is specified in various ways. We write $U_+$ for the space of non-negative excursions:
\begin{equation}
U_+ = \{ f: \R^+ \rightarrow \R^+ | f^{-1}( (0,\infty) ) = (0,\zeta)\; \hbox{\rm for some } \zeta >0 \}.
\end{equation}
We let $n_+$ be the law of excursions away from zero of $|W|$, so that 
\begin{equation}
n_+ (\{ f \in U : \sup_t f(t) > x\}) = 1/x.
\end{equation}
Now let $\Pi$ be  a Poisson random measure on $(0,\infty) \times U_+\times \{1, \ldots, k\}$ with mean measure 
\begin{equation}
M_N \equiv \kk^{-1} dt \times n_+(df) \times \mu_\kk,
\label{MNdef}
\end{equation}
  where $\mu_\kk(\{j\})=1$ for each $j=1, \ldots , N$. Now we define the clock
\begin{equation}
T(\ell) =  \iiint_{ (0,\ell] \times U_+\times S_\kk} \zeta(f) \; \Pi(ds,df, dj),
\end{equation}
with inverse
\begin{equation}
L_t = \inf\{ \ell: T(\ell) > t \}.
\end{equation}
The local time $L$ remains constant through all excursion intervals; let
\begin{equation}
 g_t \equiv \sup\{ s: L_s < L_t \}
\end{equation}
denote the left end of the excursion including time $t$. Finally, we may define
\begin{eqnarray}
X_j(t) &=& -\kk^{-1} L_t + f(t - g_t) \quad \hbox{\rm if $(L_t,f,j)$ is a point of $\Pi$;}
\label{54}
\\
&=& -\kk^{-1} L_t  \qquad \hbox{\rm else.}
\label{55}
\end{eqnarray}

\end{proof}

\begin{proposition}\label{prop4}
Uniqueness in law holds for \eqref{SDE0_bis}.
\end{proposition}
\begin{proof}
Firstly, let us deal with the case $N=2$. We saw in Section \ref{S2} that any solution to the SDE \eqref{SDE0} can be represented in terms of the Brownian motion $Y$ defined at \eqref{dY_tanaka} by the equations \eqref{X1fromY}, \eqref{X2fromY}, so the law of the solution $(X_1,X_2)$ is uniquely determined.

\medskip

The case $N\geq 3$ requires a little more subtlety. Take any $j \neq k$ in $\{ 1, \ldots, \kk\}$, and define
\begin{eqnarray}
A_t &=&\int_0^t (I_j(X(s) ) + I_k(X(s) )) \; ds,
\\
\tau_t &=& \inf\{ u: A_u > t\}.
\end{eqnarray}
Notice that
\begin{equation}
(1 - I_j(X(t))-I_k(X(t))) \, dX_j(t) = 
I_j(X(t)) (1 - I_j(X(t))-I_k(X(t))) \, dW_t = 0
\end{equation}
so $X_j$, $X_k$ do not change when the clock $A$ is not growing. Therefore if we define
\begin{equation}
\tilde{X}_j(t) = X_j(\tau_t), \quad
\tilde{X}_k(t) = X_k(\tau_t),
\end{equation}
we have
\begin{equation}
\inf\{ \tilde{X}_j(s) : s \leq t\} = \uX_j(\tau_t) = \uX(\tau_t).
\end{equation}
With a slight overloading of notation, we have
\begin{eqnarray*}
d\tilde{X}_j(t) &=& I_j(\tilde{X}(t)) \; d\tilde W_t,
\\
d\tilde{X}_k(t) &=& I_k(\tilde{X}(t)) \; d\tilde W_t,
\end{eqnarray*}
so {\em the pair $(\tilde{X}_j,\tilde{X}_k )$ is a solution of the SDE for the case $d=2$}. But we know that uniqueness in law holds for this situation, so in particular we know that at any time $t$ such that $\tilde{X}_j(t) > \tilde{X}_k(t)$ we must have
\begin{equation}
\tilde{X}_k(t)  = \inf\{ \tilde{X}_k(s) : s \leq t\} = 
\inf\{ \tilde{X}_j(s) : s \leq t\} = \uX(\tau_t).
\end{equation}
But the choice of the pair $j,k$ was arbitrary, so we deduce that
\begin{center}
whenever $X_k(t) > \uX(t)$, it must be that $X_j(t) = \uX(t)$ for all $j \neq k$.
\end{center}
{ \em This is statement (3) of Theorem \ref{thm1}.  }
In view of the facts that $\sum_j X_j = W$ and $\uX = N^{-1} \uW$, if at time $t$ we have  $X_k(t) > \uX(t)$ then
\begin{equation}
X_k(t) - \uX(t) = \sum_j \{ X_j(t) - \uX(t) \} = W_t - \uW_t,
\label{62}
\end{equation}
{\em proving statement (4) of Theorem \ref{thm1}.}

The picture of the solution of \eqref{SDE0_bis} is now clearer: {\em each excursion away from 0 of $W - \uW$ is assigned to exactly one of the $X_j$}. We shall now prove that  the probabilistic structure  of {\em any} solution to \eqref{SDE0_bis} coincides with the probabilistic structure of the solution constructed in Proposition \eqref{prop3}.

For this, define
\begin{equation}
A_j(t) \equiv \int_0^t I_j(X(s))\; ds, 
\qquad \tau_j(t) \equiv \inf\{ s: A_j(s) > t \}, \qquad X^j(t) = X_j(\tau_j(t)).
\label{dAjdef}
\end{equation}
Thus each $X^j$ is a standard Brownian motion. In fact, the $X^j$ are {\em independent} Brownian motions, as we see by the following argument. Fix any $t_1, \ldots, t_N>0$, and any $\theta_1, \ldots, \theta_N \in \R$.  Then
\begin{equation*}
M_t \equiv \exp\bigl[\;
\sum_{j=1}^N\{  i \theta_j X_j(t \wedge \tau_j(t_j))-\half 
\theta_j^2 A_j(t \wedge \tau_j(t_j))  \}
\;\bigr]
\end{equation*}
is a bounded martingale\footnote{This is proved using It\^o's formula, and the fact that $d \left\langle X_j, X_k \right\rangle =0$ for $j \neq k$.}, so
\begin{equation*}
1 = EM_0 = EM_\infty = 
\exp\bigl[\;
\sum_{j=1}^N\{  i \theta_j X^j(t_j)-\half 
\theta_j^2 t_j  \}
\;\bigr].
\end{equation*}
Hence the $X^j(t_j)$ are independent zero-mean Gaussian, and the independence of the $X^j$ follows. Thus if we decompose each $X^j$ into its Poisson process $\Pi^j$ of excursions away from the minimum as in Proposition \ref{prop3}, then the $\Pi^j$ are independent. Therefore if we define a random measure $\bar{\Pi}$ on $(0,\infty) \times U_+\times \{1, \ldots, k\}$ by
\begin{equation*}
 \bar{\Pi}(B \times \{j\} ) = \Pi^j(B)
 \end{equation*}
for any Borel $B \subseteq (0,\infty) \times U_+$, it can be seen that $\bar{\Pi} $ is a Poisson random measure, with mean measure $ dt \times n_+(df) \times \mu_\kk$, where as before $\mu_N(\{j\})=1$ for $j=1, \ldots,N$. This is the measure $M_N$ defined at \eqref{MNdef}, but scaled up by a factor of $N$. A point in $(0,\infty) \times U_+\times \{1, \ldots, k\}$ is a triple, where the final component in $\{1, \ldots,N\}$ we refer to as the {\em label}. If we take all points in $\bar{\Pi}$ with label $j$, we see the Poisson point process of excursions of a Brownian motion (in fact, $X^j$); if we take all points in $\bar{\Pi}$, we see the Poisson point process of a Brownian motion (in fact, $W$) but scaled up by a factor of $N$. This means that the local time of the totality of all points in $\bar\Pi$ is growing $N$ times as fast as the local time of the corresponding Brownian motion, explaining the factor $N^{-1}$ in expressions \eqref{54}, \eqref{55}.

\end{proof}

\begin{proposition}\label{prop5}
For any $k$, the process
\begin{equation}
M_t \equiv 
X_k(t) - (N-1)^{-1} \sum_{j \neq k} X_j(t)= \frac{NI_k(X(t)) -1}{N-1}  (W(t) - \uW(t) )
\label{mart5}
\end{equation}
is a martingale.
\end{proposition}
\begin{proof}
Firstly we verify the algebraic equivalence of the two sides of \eqref{mart5}. If $X_k$ is the lead process at time $t$ (that is,  $I_k(X(t))=1$), then from \eqref{62} we have $X_k(t) = W(t) - \uW_t + \uX(t)$, and $X_j(t) = \uX(t)$ for $j \neq k$, so the  two sides of \eqref{mart5} agree in this case. If $X_k$ is not the lead process at time $t$, then similarly the left-hand side of \eqref{mart5} is equal to 
$ - (N-1)^{-1}(W(t) - \uW(t))$, as required.

Now take $0 < s < t$, any $A \in \sF_s$, and let $\tau = \inf\{u>s: W(u) = \uW(u)\}$.
Notice that
\begin{eqnarray*}
(N-1)E[M_t - M_\tau: A, \tau<t] &=& E[ (NI_k(X(t))-1)(W(t) - \uW(t)): A, \tau<t]
\\
&=& E[ NI_k(X(t))-1:A, \tau<t]\; E[W(t) - \uW(t): A, \tau<t]
\\
&=& 0
\end{eqnarray*}
since the label of any excursion is independent of the path of that excursion, and each label has equal probability $1/N$.
  Therefore
\begin{eqnarray*}
E[ M_t - M_s : A] &=& E[M_t-M_s:A, \tau \leq t] + E[M_t-M_s:A, \tau > t]
\\
&=& E[M_\tau-M_s:A, \tau \leq t] + E[M_t-M_s:A, \tau > t]
\\
&=& E[M_{\tau\wedge t}-M_s:A]
\\
&=& E[ (N-1)^{-1}(I_k(X(s))-1) (W_{\tau\wedge t}-W_s) : A]
\\
&=& E[  E(W_{\tau\wedge t}-W_s | \sF_s)\; (N-1)^{-1}(I_k(X(s))-1)I_A]
\\
&=& 0,
\end{eqnarray*}
because the label of the lead process does not change during $[s,\tau]$, nor does $\uW$.

\end{proof}

This completes the proof of Theorem \ref{thm1} in the case where $X(0)=0$. The general case follows by concatentation. So if we have $X_1(0) > X_2(0) > \ldots >X_N(0)$, then up until the time $T_j \equiv \inf\{ t: X_1(t) = X_j(0)\}$ none of the processes $X_i, \; i \geq j$ has moved. Up until $T_2$, only $X_1$ is moving, so this behaves like Brownian motion. Between $T_2$ and $T_3$, both $X_1$ and $X_2$ are moving with a common minimum, so we may apply Theorem \ref{thm1} with two processes, both starting at the same place; then between $T_3$ and $T_4$ we have three moving processes, and so on.

\end{document}